\documentclass[11pt]{amsart}
\theoremstyle{plain}
\newtheorem{Thm}{Theorem}

\begin{document}

\title[Uniqueness of Schrodinger flow]
{Uniqueness of Schrodinger flow via energy inequality}

\author{Li MA, Lin Zhao, and Jing Wang}

\address{Department of mathematical sciences \\
Tsinghua university \\
Beijing 100084 \\
China} \email{lma@math.tsinghua.edu.cn} \dedicatory{}
\date{Feb. 20th, 2008}

\begin{abstract}

In this short note, we show a uniqueness result of the energy
solutions for the Cauchy problem of Schrodinger flow in the whole
space $R^n$ provided there is a smooth solution in the energy class.


\end{abstract}

\thanks{$^*$ The research is partially supported by the National Natural Science
Foundation of China 10631020 and SRFDP 20060003002 }
 \maketitle

\section{Introduction}

Uniqueness problem for weak solutions of nonlinear evolution system
is an important and hard topic. In this short note, we study the
uniqueness result for Schrodinger (map) flow from
$\mathbf{R}^2\times [0,+\infty)\to S^2$ (since we have no intention
to make survey on this topic, we refer to
\cite{Landau},\cite{Ding88}, and \cite{BIKT08} for more physical
background and results from K.Uhlenbeck, C.Terng, J.Shatah,
C.E.Kenig, T.Tao, D.Tataru, A.D.Ionescu, and others). Here $S^2$ is
equipped with the standard metric. By definition, the Cauchy problem
of the Schrodinger map flow  is a smooth mapping
$$
u: \mathbf{R}^2\times [0,+\infty)\to S^2
$$
satisfying
\begin{equation}\label{sch}
u_t=u\times \Delta u, \; \; in \; R^2\times [0,+\infty) ,
\end{equation}
with the Cauchy data $$ u|_{t=0}=u_0,
$$
where $\Delta u=\sum\frac{\partial^2 u}{\partial x_i^2}$ is the
usual Laplacian operator in $R^2$. Here and thereafter, we use the
sum convention. From the equation (\ref{sch}) we can easily see that
for the regular solution $u$, the energy $$
E(u(t))=\frac{1}{2}\int_{\mathbf{R}^2} |\nabla u(t)|^2dx =E(u(0))
$$
is conserved. In fact, from the fact
$$
(u_t,\Delta u)=(u\times \Delta u, \Delta u)=0,
$$
we derive that $$ \int\frac{d}{dt}E(u(t))=\int_{\mathbf{R}^2}
(\nabla u,\nabla u_t)= -\int_{\mathbf{R}^2} (\Delta u, u_t)=0.
$$
Another useful form of the equation (\ref{sch}) is
\begin{equation}\label{sch1}
-\Delta u-u|\nabla u|^2=u\times u_t.
\end{equation}
Some people like to use another form. If we compose $u$ with a
stereographic projection $\pi$ such that $z(x,t)=\pi\circ u(x,t)\in
\mathbf{R}^2$, then we have
$$
iz_t-\Delta z=\frac{2\bar{z}}{1+|z|^2}\partial_jz\partial_j z.
$$

We say that the mapping
$$v\in L^1_{loc}(\mathbf{R}^2\times [0,+\infty);
S^2)$$ is an energy class solution (and also called the weak
Schrodinger flow) if $v$ satisfies (\ref{sch}) in the distributional
sense and the energy inequality
$$
E(v(t))\leq E(v(s)), \; \; for \; all \; t>s\geq 0.
$$

We denote by $\mathbf{R}_+=[0,+\infty)$, $Du=(\nabla u,u_t)$, and
$(A,B,C)$ the determinant of the matrix formed by the ordered column
vectors $A,B,C$. We have the following uniqueness result.

\begin{Thm}\label{uniq0}
Let $$u\in C^{\infty}(\mathbf{R}^2\times [0,+\infty); S^2)$$ with
$$
|Du|\in L^{\infty}_{loc}([0,\infty), L^{\infty}(\mathbf{R}^2))
$$
and
$$
v\in L^1_{loc}(\mathbf{R}^2\times [0,+\infty); S^2), \;  with \;
 |Dv|\in L^{\infty}_{loc}([0,+\infty); L^{\infty}(\mathbf{R}^2))
$$
 be (weak) solutions to (\ref{sch}) with the smooth Cauchy data
$u_0$.

Assume that $|\nabla u_0|$ has compact support. Then $u=v$.
Furthermore, the same result is true for any Euclidean space
$\mathbf{R}^n$ in place of the plane $\mathbf{R}^2$.
\end{Thm}
 We remark that the assumption about the $L^{\infty}$
bound of the t-derivative $v_t$ is equivalent to $|\Delta v|(t)\in
L^{\infty}(\mathbf{R}^2)$ by using the Schrodinger map flow equation
(\ref{sch1}). We remark that Theorem \ref{uniq0} in the plane case
can also be derived from a result in \cite{kenig05}.  It seems to us
that $L^{\infty}$ bound of the t-derivative $v_t$ is very stronger.
However, our assumption in dimension two can be weaken by using the
Ladyzhenskaya's inequality. Recall here that the Ladyzhenskaya's
inequality (\cite{Lady}) is the following one. For any $w\in
C^1_0(\mathbf{R}^2)$,
 $$ ||w|^2|_{L^2(\mathbf{R}^2)}\leq
C|w|_{L^2(\mathbf{R}^2)}|\nabla w|_{L^2(\mathbf{R}^2)}.
$$

We point out that similar results for wave equation and wave maps
have been obtained by M.Struwe \cite{St99} and here we use the
similar Gronwall inequality method. One may also see the paper
\cite{St99} for more related results. Similar result for Schrodinger
system has been obtained in \cite{ma}. The new ingredient in our
proof is the simple $L^2$ growth of $v-u$ and its derivatives. In
the following, we denote by $C$ the various constants, which do not
depend on time variable $t$.

\section{Proof of Theorem \ref{uniq0}}
For simple presentation, we just give the proof for $n=2$. The
general case can be treated in the same way.

 Write $v=u+w$. Note that
$|u|^2=1=|v|^2$. So, $u_t\cdot u=0=v\cdot v_t$.

Consider
$$
\frac{1}{2}\frac{d}{dt}\int_{\mathbf{R}^2}|w|^2dx=\int_{\mathbf{R}^2}w\cdot
w_t.
$$
Note that
$$
w\cdot w_t=w\cdot v\times \Delta v-w\cdot u\times \Delta u=w\cdot
u\times \Delta w.
$$
Then
$$
\frac{1}{2}\frac{d}{dt}\int_{\mathbf{R}^2}|w|^2dx=\int_{\mathbf{R}^2}(w\times
u)\cdot \Delta w.
$$
Using the integrating by part, the latter term can be bounded by
$$
\leq C\int_{\mathbf{R}^2}(|\nabla w|^2+|w|^2).
$$
Hence, we have
\begin{equation}\label{L2}
\int_{\mathbf{R}^2}|w|^2dx(t)\leq
C\int_0^td\tau\int_{\mathbf{R}^2}(|\nabla w|^2+|w|^2).
\end{equation}

We clearly have that
$$
E(v)=E(u)+I+E(w),
$$
where $$ I=\langle dE(u),w\rangle=\int_{\mathbf{R}^2} \nabla u\cdot
\nabla w dx.$$ Hence,
\begin{equation}\label{key}
0\geq E(v(t))-E(v(0))=I(t)-I(0)+E(w(t))-E(w(0)).
\end{equation}

 Consider
\begin{equation} \label{III}
I(t)-I(0)=\int_0^t d\tau\int_{\mathbf{R}^2}\partial_t(\nabla u\cdot
\nabla w) dx.
\end{equation}
Then we have
$$
(\ref{III})=-\int_0^t d\tau\int_{\mathbf{R}^2}[(\Delta u)\cdot
w_t+(\Delta w)\cdot u_t],
$$
which can be written as
\begin{eqnarray*}
(\ref{III})&=&\int_0^t d\tau\int_{\mathbf{R}^2}w_t\cdot[u|\nabla
u|^2-v|\nabla v|^2+u\times u_t-v\times v_t]\\
&=&\int_0^t d\tau\int_{\mathbf{R}^2}w_t\cdot[u|\nabla u|^2-v|\nabla
v|^2-w\times u_t]\\
&=&\int_0^t d\tau\int_{\mathbf{R}^2}[w_t\cdot(u|\nabla u|^2-v|\nabla
v|^2)+(w,w_t,u_t)]
\end{eqnarray*}
after using the relations
$$
w_t\cdot\Delta u=-w_t\cdot (u|\nabla u|^2+u\times u_t),
$$
and
\begin{eqnarray*}
u_t\cdot \Delta w&=&u_t\cdot \Delta v-u_t\cdot \Delta u\\
&=& -u_t\cdot (v|\nabla v|^2+v\times v_t)\\
&=& w_t\cdot (v|\nabla v|^2+v\times v_t).
\end{eqnarray*}

We now bound each term case by case. Using $$ w_t=v\times \Delta
w+w\times \Delta u,
$$
we have
\begin{eqnarray*}
\int_0^t d\tau\int_{\mathbf{R}^2}(w,w_t,u_t)&=&\int_0^t
d\tau\int_{\mathbf{R}^2}(w,v\times \Delta w,u_t)\\
&=&-\int_0^t d\tau\int_{\mathbf{R}^2}(w\times u_t,v,\Delta w)
\end{eqnarray*}

Upon integrating by part, we have
$$
|\int_0^t d\tau\int_{\mathbf{R}^2}(w\times u_t,v,\Delta w)|\leq
C\int_0^t d\tau\int_{\mathbf{R}^2}(|\nabla w|^2+|w|^2)dx.
$$

To bound other terms, we notice the relation that $$ v|\nabla
v|^2-u|\nabla u|^2=w|\nabla v|^2+u(|\nabla v|^2-|\nabla u|^2).
$$
Using $$
 w_t=v\times \Delta
w+w\times \Delta u,
$$ again, we have
\begin{eqnarray*}
w_t\cdot(v|\nabla v|^2-u|\nabla u|^2)&=&(v\times \Delta w+w\times
\Delta u)\cdot w|\nabla v|^2+\\
&& (v\times \Delta w+w\times \Delta u)\cdot u(|\nabla v|^2-|\nabla
u|^2)\\
&=&(v\times \Delta w+w\times
\Delta u)\cdot w|\nabla v|^2+\\
&& v\times \Delta w\cdot u(|\nabla v|^2-|\nabla u|^2)\\
&&+w\times \Delta u\cdot u(\nabla w\cdot\nabla (u+v))\\
&=&(v\times \Delta w+w\times
\Delta u)\cdot w|\nabla v|^2+\\
&&w\times \Delta w\cdot u(|\nabla v|^2-|\nabla u|^2)\\
&&+w\times \Delta u\cdot u(\nabla w\cdot\nabla (u+v)) ,
\end{eqnarray*}
which implies that upon integrating by part, based on a
regularization approach,
$$
|\int_{\mathbf{R}^2}w_t\cdot(v|\nabla v|^2-u|\nabla u|^2)|\leq
C(E(w)+|w|^2_2).
$$

 Then using (\ref{key}), we have
$$
|I(t)-I(0)|\leq C\int_0^t(E(w)+|w|^2_2),
$$
and we have
$$
E(w(t))\leq E(w(0))+C\int_0^t(E(w)+|w|^2).
$$
This estimate gives us the stability result for the Schrodinger map
flow.

Note that $E(w(0))=0$. Combining this with (\ref{L2}), we have
$$
E(w(t))+|w|^2_2(t)\leq C\int_0^t(E(w)+|w|_2^2).
$$
 By the Gronwall inequality, we have
$$
w=0, \; \; and \;  u=v.
$$
This completes the proof of Theorem \ref{uniq0}.

\end{document}